\documentclass[5p,times]{elsarticle}

\usepackage{lineno,hyperref}
\modulolinenumbers[5]

\journal{Journal of Parallel Computing}

\usepackage{amsmath}
\usepackage{amssymb}
\usepackage{todonotes}
\usepackage[utf8]{inputenc}

\usepackage[absolute,overlay]{textpos}
\usepackage{setspace}

\usepackage{cuted}
\usepackage{nccmath}

\usepackage[super]{nth}

\begin{document}

\begin{frontmatter}

\title{Exponential Integrators with Parallel-in-Time Rational Approximations for the Shallow-Water Equations on the Rotating Sphere}

\author{Martin Schreiber\corref{mycorrespondingauthor}}
\ead{martin.schreiber@tum.de}
\address{
	Chair of Computer Architecture and Parallel Systems, Department of Informatics, Technical University of Munich, Germany
	\\
	Department of Computer Science / Mathematics, University of Exeter, EX4 4QF, United Kingdom
}

\author{Nathanaël Schaeffer\corref{}}
\ead{nathanael.schaeffer@univ-grenoble-alpes.fr}
\address{Univ.\,Grenoble Alpes, CNRS, ISTerre, F-38000 Grenoble}

\author{Richard Loft}
\ead{loft@ucar.edu}
\address{CISL/NCAR, Boulder, USA}

\cortext[mycorrespondingauthor]{Corresponding author}

\begin{abstract}

High-performance computing trends towards many-core systems are expected to continue over the next decade.
As a result, parallel-in-time methods, mathematical formulations which exploit additional degrees of parallelism in the time dimension, have gained increasing interest in recent years.
In this work we study a massively parallel rational approximation of exponential integrators (REXI).
This method replaces a time integration of stiff linear oscillatory and diffusive systems
by the sum of the solutions of many decoupled systems, which can be solved in parallel. Previous numerical studies showed that this reformulation allows taking arbitrarily long time steps for the linear oscillatory parts.

The present work studies the non-linear shallow-water equations on the rotating sphere, a simplified system of equations used to study properties of space and time discretization methods in the context of atmospheric simulations.
After introducing time integrators, we first compare the time step sizes to the errors in the simulation, discussing pros and cons of different formulations of REXI.
Here, REXI already shows superior properties compared to explicit and implicit time stepping methods.
Additionally, we present wallclock-time-to-error results revealing the sweet spots of REXI obtaining either an over $6 \times$ higher accuracy within the same time frame or an about $3 \times$ reduced time-to-solution for a similar error threshold.
Our results motivate further explorations of REXI for operational weather/climate systems.

\end{abstract}


\end{frontmatter}
%

\sloppy

\section{Introduction}

Solving nonlinear partial differential equations (PDEs) has long been one of the main applications of high-performance computing (HPC).
During the four decades beginning in the 1960s, improvement in computational performance of HPC systems was primarily driven by the increase of processor clock frequency.
Sometimes called the ``free lunch''\cite{sutter2005free} era, new hardware could be relied on to provide new capabilities.
In particular, these new capabilities directly led to improvements in the time-to-solution without requiring any algorithmic changes.

However, since the mid 2000s, solving PDEs efficiently faces a new challenge:
the stagnation and even decline of processor clock frequency, and a compensatory increase in processor parallelism.
As a consequence, a fundamental redesign is required to exploit new ways to advance the quality of simulations.

Our work is motivated by PDE solvers for climate and weather simulations which have long been considered grand challenge problems in HPC.
In the present work, we put the focus on the parallelization of a rational approximation of exponential integrators (REXI).
This method allows a decomposition of the computations of one time step (which are typically strictly sequential computations), into a set of terms which can be solved independently, followed by an accumulation of the results.
We study this new method in comparison with various other standard time-stepping methods, focusing on how the errors and the wallclock time behave once utilizing the potential to parallelize over the independent terms of REXI.
In this way we can clearly see how parallel-in-time methods lead to higher accuracy by exploiting additional degrees of parallelization as part of the time integration.
To the best of our knowledge, this is the first time that such studies are conducted in parallel on a supercomputer for the full non-linear shallow-water equations on the rotating sphere using a formulation based on a massively parallel rational approximation.

\section{Related work}
\label{sec:related_work}

We first provide a brief overview of exponential integrators (Sec.~\ref{sec:exp_integrators}) followed by parallel-in-time methods (Sec.~\ref{sec:pint_methods}) and the relation to exponential integrators.
With the focus on climate and weather applications, we close the related work section with an overview of state-of-the-art time-stepping methods in the field (Sec.~\ref{sec:timestepping_methods_for_c_and_w}).

\subsection{Exponential integrators}
\label{sec:exp_integrators}

Exponential integrators (EI) have been studied for several decades\cite{Moler2003,hochbruck2010exponential}:
for generic introductions into EIs and for various formulations of EIs, see \cite{lawson1967generalized,hochbruck1998exponential,CoxMathews2002,tokman2006efficient,tokman2011new}.
Most relevant to the present paper, Clancy et al.~\cite{clancy2011laplace} applied EI in the context of Laplace transforms to gain higher accuracy. Spherical harmonics were exploited to solve for each term in the Laplace operator formulation as part of a REXI formulation.
EIs have also been applied with Krylov subspace solvers (see \cite{Niesen2012} for an overview).
Clancy et al.~\cite{ClancyPudykiewcz2013} used such a solver strategy, but lost the possibility to exploit additional parallelism due to the purely sequential iterations over the Krylov subspaces.
Bonaventura \cite{Bonaventura15} studied local EIs by exploiting space-time locality properties of physical phenomena.
There, additional degrees of freedom were generated by decomposing the originally globally connected domain into multiple independent, localized overlapping subdomains.
Garcia et al.~\cite{garcia2014exponential} studied EIs for thermal convection in rotating spherical shells using Krylov approximations.
They assessed the performance of their algorithm in terms of the ratio of error to time-step size for thermal convection.

\subsection{Parallel-in-time methods}
\label{sec:pint_methods}

Parallel-in-time (PinT) methods seek to exploit an increasing parallelism of computational resources, by allowing computation of different time steps in the overall time integration in parallel.
An overview of 50 years of research in this direction can be found in the work of Gander \cite{Gander15} and, for a brief review, we only focus on the PinT developments which show the most significant relevance to our work.
The use of computing resources in this way is appropriate for throughput-driven problems where not meeting time-to-solution or error thresholds can make results less valuable (e.g.\,augmented surgeries) or in the worst case even entirely useless (e.g.\,weather forecast).
Hence, parallel-in-time methods should be used only when the spatial scalability is saturated.
The Parareal algorithm formed the basis of current research on PinT methods and uses an iterative-in-time scheme \cite{Lions2001}.
A related method is the Parallel Full Approximation Scheme in Space and Time (PFASST) method which is based on spectral deferred corrections in time which successively increases the integration order in combination with a multi-resolution approach\cite{emmett2012toward}.
In this work we focus on direct parallel methods.
One of such methods is based on revisionist integral deferred correction method (RDIC)\cite{christlieb2010parallel}
which consists of stages executed in parallel with built-in corrections of previous results.
However, explicit RDIC time stepping suffers of a reduced stability region compared to higher-order explicit Runge-Kutta based time steppers.

EI formulations in the context of parallel-in-time were first applied with Paraexp (see Gander et al.\,\cite{Gander2013}).
EIs also represent one of the main building blocks in the asymptotic PinT approach (Haut et al.\,\cite{Haut_Wingate_14} who used an analytical solution of the EI for the linear parts of the shallow-water equations on the plane).
All these formulations require solving exponentials of matrices which can be accomplished with a Rational Approximation of Exponential Integrators (REXI) and its variants \cite{clancy2011laplace,ClancyPudykiewcz2013,Haut2015}.
For exponentials of linear operators, such rational approximations can be used to approximate time-integrating functions with a sum over solutions for a complex-valued system of linear equations.
Each solution can be computed totally independently, hence parallel in time, and can be used as fundamental building blocks for exponential integrators.
In particular, for linear operators, REXI allows arbitrarily long time steps and precision beyond the timestepping accuracy of standard methods, which served as an initial motivation of the present work.
Haut et al.\,also developed a computationally efficient and stable method \cite{Haut2015} to compute the weights $\alpha_n$ and $\beta_n$ to solve the EI formulation for purely oscillatory problems (e.g.\,the linearized shallow-water equations) and shares its algebraic structure $\sum_n \beta_n(\Delta t L+\alpha_n)^{-1} U(0)$ with the previously mentioned Laplace transforms \cite{clancy2011laplace}.
Parallel performance and numerical studies of this method for the shallow-water equations on the plane were conducted by Schreiber et al.\,\cite{schreiber2016beyond} using spectral solvers with the Fourier space and finite-difference formulations.
Significant speedups of two orders of magnitude for linear problems were obtained there by exploiting additional degrees with a rational approximation of oscillatory problems and a spectral fast Helmholtz solver.
Various ways to incorporate non-linear parts in exponential integrator formulations exist (see \cite{lawson1967generalized,tokman2006efficient,tokman2011new,CoxMathews2002}). In this work, we use Strang-splitting and EI formulations denoted Exponential Time Differencing $n$-th order Runge Kutta (ETDnRK) \cite{CoxMathews2002} where we use a REXI parallel-in-time approach to solve the exponentials of matrices.

\subsection{Time stepping methods for climate and weather simulations}
\label{sec:timestepping_methods_for_c_and_w}

We first give a brief overview of the two most significant time integration schemes employed in the dynamical cores\footnote{Here we interpret a dynamical core as the software component of climate/weather simulations for grid-resolved processes, without external physics effects.} for weather and climate simulations, which will be compared to (parallel-in-time) EIs.
The first of these are splitting methods, which treat various equation terms differently based on their mathematical and physical properties.
One of the first big advancements was to use implicit methods for the fast gravity modes.
This was also exploited by Robert \cite{robert1969integration} using special properties of spherical harmonics.
Implicit time integration led to stable time-step sizes which are $7\times$ to $10\times$ larger if compared to the maximum stable time step size for explicit time integration methods.
However, larger time step sizes than this were not possible due to stability limitations with the explicitly treated nonlinearities.
Robert \cite{robert1982semi} applied semi-Lagrangian methods in combination with semi-implicit time stepping which led to another increase in time step sizes of $4 \times$ -- $6 \times$.
Such semi-Lagrangian methods are currently also used in operational weather simulation codes, see e.g.\,\cite{Barros1995,wood2014inherently}.
In the present work, we will put our focus solely on semi-implicit methods.

\section{Shallow-water equations on the rotating sphere}

\label{sec:shallow_water_equations}

Studying new discretization methods for global climate and weather simulations is frequently first tested with the shallow-water equations (SWE), hence ignoring the vertical discretization and reducing the original problem to a single-layer simulation (see e.g.\,\cite{williamson1992standard} for standardized tests).
The SWE represent important physical properties of a single layer of the full atmospheric model to assess the numerics.
For the studies in this work, this allows us to focus purely on time integration of the horizontal parts and we start by briefly introducing the SWE on a rotating sphere.
We use the vorticity-divergence formulation (see \cite{temperton1991scalar,hack1992description})
\begin{eqnarray}
	\frac{\partial\Phi'}{\partial t}	&=&  -\overline{\Phi}\delta 	-\nabla  \cdot (\Phi' \mathbf{V})\\
	\frac{\partial\zeta}{\partial t}	&=&	-f\delta-\mathbf{V} \cdot (\nabla f) -\nabla \cdot (\zeta\mathbf{V})\\
	\nonumber
	\frac{\partial\delta}{\partial t}	&=&  f\zeta+\mathbf{k} \cdot (\nabla f)\times\mathbf{V}-\nabla^{2}\Phi\\
										& &  + \mathbf{k} \cdot \nabla\times(\zeta\mathbf{V})-\nabla^{2}\frac{V \cdot V}{2}
\end{eqnarray}
using the vorticity $\zeta$, divergence $\delta$ and the perturbation of the geopotential $\Phi' = h' g$ with $h'$ the perturbed surface height, $\overline{\Phi}$ the global average of the geopotential, $g$ the acceleration of gravity, $f=2\Omega \sin \phi$ the Coriolis effect varying along the latitude $\phi$ and $\textbf{V} = (u,v)^T$ the velocity.
Furthermore, we use the identities $\mathbf{V}=\mathbf{k}\times\nabla\psi+\nabla\chi$ with the stream $\psi = \nabla^{\text{-2}}\zeta$ and potential $\chi=\nabla^{\text{-2}}\delta$ to compute the velocities from the vorticity and divergence.
Given the velocity, we can compute the divergence $\delta=\nabla \cdot \mathbf{V}$ and vorticity $\zeta=\mathbf{k} \cdot (\nabla\times\mathbf{V})$ with $\mathbf{k}$ the unit vector orthogonal to the surface.
Such transformations between velocity and stream/vorticity formulation avoid inconsistencies of vector fields close to the poles if transformed to/from spectral space with scalar spherical harmonics.

Using splitting methods to integrate parts of the SWE differently, we start with a nomenclature of the different terms in the SWE.
We split up the right-hand side into its subcomponents ($L_g$, $L_c$, $N$) with $U$ representing all state variables ($\Phi'$, $\zeta$, $\delta$):

\begin{strip}
\begin{ceqn}
\begin{eqnarray}
\underbrace{
\begin{array}{c}
	\frac{\partial\Phi'}{\partial t}\\
	\frac{\partial\zeta}{\partial t}\\
	\frac{\partial\delta}{\partial t}
\end{array}
}_{\frac{\partial U}{\partial t}}
\begin{array}{c}=\\=\\=\end{array}
\underbrace{
\begin{array}{c}
-\overline{\Phi}\delta\\
0\\
-\nabla^{2}\Phi
\end{array}
}_{L_g U}
\underbrace{
\begin{array}{c}
+0\\
-f\delta-\mathbf{V} \cdot (\nabla f)\\
+f\zeta+\mathbf{k} \cdot (\nabla f)\times\mathbf{V}
\end{array}
}_{L_c U}
\underbrace{
\begin{array}{c}
-\mathbf{V} \cdot \nabla \Phi'-\Phi' \nabla \cdot  \mathbf{V}\\
-\nabla \cdot (\zeta\mathbf{V})\\
+\mathbf{k} \cdot \nabla\times(\zeta\mathbf{V})-\nabla^{2}\frac{\mathbf{V} \cdot  \mathbf{V}}{2}
\end{array}.
}_{N(U)}
\end{eqnarray}
\end{ceqn}
\end{strip}

For the realization of the different time splitting approaches in the following sections and a compact notation, we write
\begin{equation}
	\frac{\partial U}{\partial t} = L_g U + L_c U + N(U).		\label{eq:sw_compact}
\end{equation}
All of these terms hold different mathematical and physical properties (see e.g.\,\cite{lemarie2015stability}) and we briefly describe their physical relevance.

\subsection{Physical relevance}
For the linear parts, the gravity term $L_g$ relates to hyperbolic wave propagations which are generated due to gravitational force and the resulting linear propagation which is the most restrictive term for time step sizes using explicit time integrators.
The (linear) Coriolis effect $L_c$ is generated by the rotation of the Earth, generating amongst others important balancing effects in the atmosphere, the so-called geostrophic balance.
Since this effect also plays a crucial role in the propagation of Rossby waves, we expect also improved accuracy by accurately tracking them.
The non-linear term $N(U)$ includes also the non-linear advection term which will be the main limiting factor of all time stepping methods discussed in the present work if using non-explicit time integrators for the linear parts.

\subsection{Space discretization with Spherical Harmonics}
\label{sec:spherical_harmonics}

We solve the non-linear SWE on the rotating sphere using a pseudo-spectral method on the sphere relying on spherical harmonics (SH).
The non-linear terms are computed in physical space, and thus require transforming back and forth from spectral to physical space.
The main motivation for using a spectral method is its high accuracy with respect to the non-truncated modes.
We can thus safely assume that the errors are mostly due to the time-stepping schemes that we want to study.

In this work, all transforms are computed using the efficient library SHTns \cite{schaeffer2013efficient}.
This library has been highly optimized and designed for numerical simulations of fluids in spheres \cite{gastine2015,schaeffer2017}, but is used across a variety of fields like atmospheric and climate sciences \cite{augier2013,dawson2016windspharm}, astrophysics \cite{rincon2017supergranulation} or acoustics \cite{carley2016fast}.
It has also been employed to solve the SWE \cite{suhas2017tropical}.
All SH transformations support threaded parallelization which is used throughout the present work.

Regarding the required transformations in REXI, complex-valued spatial fields are required.
The complex-valued and threaded parallelization for such transformations was developed as part of this work and we study their performance in Sec.\,\ref{sec:parallelization_in_space}.

\subsection{Exponential time integration for SWE}
	Exponential integrators allow for solving the linearized SWE ($L_g$) on the rotating sphere (including $L_c$) with very high accuracy \cite{clancy2011laplace,Schreiber2017SHREXI}.
	Furthermore, exponential integrators have been shown to support time steps as long as 1.5 days \cite{Schreiber2017SHREXI}.
	It has also been shown that SH provides an ideal basis for REXI to time integrate the linear parts of the SWE\cite{clancy2011laplace}, even on the rotating sphere \cite{Schreiber2017SHREXI}.
	However, the non-linearities $N(U)$ in the SWE pose additional restrictions regarding the stable time integration size, the resulting accuracy and finally the wallclock time which will be studied in the present work.
	Therefore, the main focus of the present work is on assessing the influence of the non-linear terms in combination with REXI.
	Details about exponential integrators will be further discussed in the following section.

\section{Time discretization}

In Section \ref{sec:shallow_water_equations}, we introduced a splitting of the different components of the SWE into multiple terms. This section describes different time splitting methods and we discuss how they are applied to each term and finally grouped together for the integration of the entire PDE.
In a general formulation, the system of equations which we integrate in time is given by
\begin{equation}
	\frac{\partial U}{\partial t}={L} U+ {N} (U)		\label{eq:split}
\end{equation}
where $U$ is the current state of the system, ${L}$ a discretized linear operator and ${N}$ denotes the terms which are assumed to be non-linear.
All our systems are considered to be autonomous, hence in particular ${L}$ is constant over time for the applicability of exponential integrators. 

\subsection{Linear single- and multi-stage integrators}

We utilize a combination of single- and multi-stage methods for the time splitting approaches and assume that the reader is familiar with explicit \nth{1} order Euler and \nth{2} order Runge-Kutta methods (see also \cite{durran2010numerical}).
The present work will be based on the Crank-Nicolson (CN) method and is frequently used in operational climate and weather simulations\footnote{Here we would like to mention that the particular formulation in each dynamical core to incorporate linear and non-linear parts can strongly vary}.
It consists of one explicit time step (forward Euler but typically with a step larger than its stability limitation), followed by an implicit one (backward Euler) which can be efficiently solved using SH.
Even if it can be computed by the combination of two single-step methods, it belongs to a special class of \nth{2} order implicit Runge-Kutta, hence multi-step, methods.
For a linear system of equations $L$ it reads
\begin{eqnarray}
	U^{n+1} =  \underbrace{\left( I - a \Delta t L \right)^{-1}}_{\text{Backward Euler}} \underbrace{\left( I + (1-a) \Delta t L \right)}_{\text{Forward Euler}} U^{n}
\end{eqnarray}
with $a = \frac{1}{2}$ throughout the paper, hence no filter.

\subsection{Rational approximation of exponential integrators (REXI)}
\label{sec:introduction_to_rexi}

We give a brief recap of the related work of rational approximations of exponential integrators for linear operators, hence it discusses time integrators for particular linear operators ($L$, $L_c$, $L_g$) from the previous section.
We start with an approximation of a function $f(x)$ which will be related to the exponential integration of an ODE.
Using a rational approximation leads to a massively parallel way to approximate $f(x)$ within an interval of approximation (related to $\alpha_{n}$) via
$
	f(x) = \sum_{n=1}^{N}\frac{\beta_{n}}{x+\alpha_{n}}
$
where $\alpha_n,\beta_n \in \mathbb{C}$ are coefficients determined later.
Here, the additional degrees of parallelization originate from the independent terms in the sum.

For linear PDE operators we can use the EI formulation
$$
	\exp(\Delta t L) = Q \exp(\Delta t \Lambda) Q^{-1}
$$
with $Q$ a matrix with the eigenvectors of $L$ and $\Lambda$ a diagonal matrix storing the associated eigenvalues.
The term $\exp(\Delta t \Lambda)$ denotes the matrix with $\exp(\Delta t \lambda_i)$ on the diagonal with $\lambda_i$ denoting the $i$-th diagonal element of the Eigenvalue matrix $\Lambda$.
Using the REXI approximation $\exp(\Delta t \Lambda) \approx 	\sum_{n=1}^{N}\beta_{n} \left( \Delta t \Lambda + \alpha_{n} \right)^{-1}$, this eventually leads to
\begin{eqnarray}
	\label{eq:rexi}
	U(\Delta t) = \exp(\Delta t L) U(0) \approx
	\sum_{n=1}^{N}\beta_{n} \left( \Delta t L + \alpha_{n} \right)^{-1} U(0).
\end{eqnarray}

The results in this work are based on REXI coefficients which will be computed using the Cauchy Contour integral method.
This method is used in the present work due to its superior properties by (i) requiring less terms compared to Haut et al. \cite{Haut2015} and (ii) being more flexible, e.g.\,directly applicable to diffusive problems.
In particular, less REXI terms led to significantly reduced resource requirements, hence a reduction in both resource consumption and parallelization overheads.
We briefly reintroduce Cauchy contour integral methods (see also e.g.\,\cite{kassam2005fourth,buvoli2015class}) to discuss the issue of cancellation errors and its mitigation developed for the present work which is required to allow large time-step sizes and to make REXI competitive to other time stepping methods.

\subsubsection{Cauchy Contour integral}
\label{sec:cauchy_contour_integral}

For a (complex-valued) holomorphic function $f(x)$ the value of the function at a particular point $x_{0} \in \mathbb{C}$ can be computed by the Cauchy contour integral
\begin{eqnarray}
	f(x_0)&=&\frac{1}{2\pi i}\ointop_{\Gamma}\frac{f(z)}{z-x_0} dz
	\label{eq:cauchy_contour_integral}
\end{eqnarray}
with $\Gamma$ the closed contour in the complex plane which encloses $x_0$.
In terms of time integration methods, this contour should enclose all eigenvalues (EVs) on the complex plane which are relevant for oscillatory (imaginary EVs) as well as diffusive (real negative EVs) linear models.

We investigate a circle as the most canonical contour which leads to a good balance between the maximum number of poles and the desired accuracy (see e.g.\,\cite{buvoli2015class}) for a relatively small radius $R < 10$.
We get the contour
\begin{eqnarray}
	\Gamma=\{R\exp(i\theta)+\mu|\theta\in[0;2\pi]\}
\end{eqnarray}
with $\mu$ a shift discussed later.
Using this in Eq.\,\ref{eq:cauchy_contour_integral} and a trapezoidal rule leads to
\begin{eqnarray}
	f(x)&=&\frac{1}{2\pi}\int_{0}^{2\pi}\frac{f(R\exp(i\theta)+\mu)\left(R\exp(i\theta)\right)}{\left(R\exp(i\theta)+\mu\right)-x} d\theta\\
		&\approx & \frac{1}{N}\sum_{n=1}^{N} \frac{f(R\exp(i\theta_n)+\mu)\left(R\exp(i\theta_n)\right)}{\left(R\exp(i\theta_n)+\mu\right)-x}.
\end{eqnarray}
Finally, we can write this compactly as a rational approximation for a function
$
	f(x) \approx \sum_{n=1}^{N}\frac{\beta_{n}}{x+\alpha_{n}}
$
with
\begin{eqnarray}
	\label{eq:cauchy_alpha}
	\alpha_n &=& -\left(R\exp(i\theta_n)+\mu\right)\\
	\label{eq:cauchy_beta}
	\beta_n  &=& - \frac{1}{N}\left(R\exp(i\theta_n)\right) f\left( R\exp(i\theta_n)+\mu \right).
\end{eqnarray}

\subsubsection{Numerical cancellation}

Using Eq.\,\ref{eq:rexi} directly results in significantly increasing numerical errors for $R \gg 10$ and a circle centered at the origin, that is $\mu=0$.
We briefly investigate this issue and a solution for physically relevant oscillatory and diffusive systems, hence excluding systems with anti-diffusion behaviour.

First, we investigate potential cancellation effects in $\alpha_n$ (see Eq.\,\eqref{eq:cauchy_alpha}) where we can directly infer $\max_{n}|\alpha_n| \approx R$, hence the magnitude of the values in $\alpha_n$ are proportional to $R$.
Second, we investigate $\beta_n$ (see Eq.\,\eqref{eq:cauchy_beta}) as the only remaining source of such errors.
For the first term, $\max_{\theta_n}|-R\exp(i\theta_n)| \approx R$ holds, hence also this term is linearly limited in its growth for larger $R$.
Next, we investigate the contribution of the function $f(x)$ to be approximated.
For $f(x) = \exp(x)$ for EI of ODEs, large real values of $x$ would lead to exponential growth.
In fact, we can observe that increasing the radius $R$ directly leads to an exponential growth of the $\beta_n$ coefficients:
$\max_{\theta_n}|Re(\exp\left(R\exp(i\theta_n)\right)))| = \exp(R)$.
We can overcome this for diffusive and oscillatory problems by a change of the contour to avoid handling of very large positive, and for our case irrelevant, eigenvalues as follows:
We first restrict the circle contour by one point $p_{0}$ with $Im(p_{0})=0$ and $Re(p_0) \geq 0$ with this point related to the cancellation errors.
The limits of the spectrum for purely oscillatory problems along the imaginary axis are required to be specified by $p_{\pm1}$ with $Re(p_{\pm1}) = 0$.
With three points on the contour of the circle, we can finally infer a radius
$r=\frac{Re(p_{0})^{2}+\left(Im(p_{-1})\right)^{2}}{2Re(p_{0})}$
and the center of the circle $p_{c} = p_{0}-R$.
For $R \gg 10$, this compensates for cancellation effects, but also increases the length of the contour of the integral.

\subsection{Non-linear exponential integrator scheme and ETDnRK}
\label{sec:etdnrk}

Exponential integrators for non-linear equations \eqref{eq:split} can be written as
\begin{eqnarray}
	U^{n+1}&=&e^{\Delta t L}U^{n}+e^{\Delta t L}\int_{0}^{\Delta t}e^{-\tau L}N(U(\tau))d\tau
\end{eqnarray}
where the dependency of the integral on $U(\tau)$ makes this non-trivial to evaluate.
ETDnRK (see \cite{CoxMathews2002}) provides one way of a higher-order approximation.
With the focus on \nth{2} order methods in this work, the time stepping is performed by 
\begin{eqnarray}
	A^{n} &=& \psi_{0}(\Delta tL)U^{n}+\Delta t\psi_{1}(\Delta tL)N(U^{n},t^{n})\\
	U^{n+1} &=& A^{n}+\Delta t\psi_{2}(\Delta tL)\left(N(A^{n},t^{n}+\Delta t)-N(U^{n},t^{n})\right)
\end{eqnarray}
with $\psi_0(K) = \exp(K)$, $\psi_1(K) = K^{-1}(\exp(K)-I)$ and $\psi_2(K) = K^{-2}(\exp(K)-I-K)$.
To infer the REXI coefficients using the Cauchy contour integral method, we use L'Hôpital's rule in the proximity of the singularities of $\psi_{1/2}(M)$ to evaluate these functions.
In a similar way, this makes REXI also applicable to all ETDnRK schemes with $n \in \{1, 2, 3, 4\}$.

\subsection{Non-linear Strang-splitting}
\label{sec:strang_splitting}

For the differential equations considered in this work
we require ways to combine different time integrators and we use
Strang-splitting\cite{strang1968construction,durran2010numerical} in the time dimension.
We split the right hand side into two separate terms
\begin{eqnarray}
	\frac{\partial U}{\partial t} = F_1(U(t)) + F_2(U(t))
\end{eqnarray}
and treat each term $F_1$ and $F_2$ individually.
In the following, we use the notation $F_a( F_b( U(t))) = F_a \circ F_b \circ U(t)$.
For a \nth{2} order accurate method, we get
\begin{equation}
	T(U^n) := U^{n+1} = F^{\frac{\Delta t}{2}}_1 \circ F^{\Delta t}_2 \circ F^{\frac{\Delta t}{2}}_1 \circ U^{n}
\end{equation}
where each term $F^{\Delta t}_n$ is integrated by an arbitrary \nth{2} order accurate time stepper with a time step size $\Delta t$.
With such a Strang-split scheme, we have the freedom to set the correspondence between $F_1$ and $F_2$ on one hand, and the linear and non-linear terms in Eq.\,\eqref{eq:split} on the other hand. This means we can choose which term is evaluated twice at each time step.
This impacts the accuracy and stability of the time stepping methods and we study this with both versions:

\textit{Version 0} integrates the linear parts twice with $F_1 = L$ and the non-linearities only once by using $F_2 = N$ with $*$ denoting a placeholder
\begin{eqnarray}
	T_{*\_ver0}(U^n) := U^{n+1} =  L^{\frac{\Delta t}{2}} \circ N^{\Delta t} \circ L^{\frac{\Delta t}{2}} \circ U^{n}.
\end{eqnarray}

\textit{Version 1} integrates the non-linear parts twice with $F_1 = N$ and the linearities only once by using $F_2 = L$
\begin{eqnarray}
	T_{*\_ver1}(U^n) := U^{n+1} = N^{\frac{\Delta t}{2}} \circ L^{\Delta t} \circ N^{\frac{\Delta t}{2}} \circ U^{n}.
\end{eqnarray}
Considering the computational complexity in terms of wallclock time and larger time step sizes for the $F_1$ term due to its splitting in two half time steps will play a significant role in the time-to-solution performance of the two versions.

\subsection{Nomenclature for time splitting methods}

\begin{table*}
	\center
	\begin{center}
		\begin{tabular}{|c|l|l|}
		\hline
			\textbf{Symbol}		&	\textbf{Description}	&	\textbf{ID}\\
		\hline
		\hline
			$L_g$	&	Linear wave motions induced by gravitational force	&	\texttt{lg}	\\
		\hline
			$L_c$	&	Coriolis effect acting on velocity components	&	\texttt{lc}	\\
		\hline
			$L = L_g + L_c$	&	Both linear terms&	\texttt{l}	\\
		\hline
			$N$		&	Non-linear advection and divergence				&	\texttt{n}\\
		\hline
		\end{tabular}
		\caption{\label{tbl:swe_terms}
			Description and string identifiers of different terms in the shallow-water equations on the rotating sphere.}
	\end{center}
\end{table*}

\begin{table*}
\center
\begin{tabular}{|l|c|}
	\hline
	Time stepping 			&Identifier\\
	\hline
	\hline
	\nth{2} order explicit Runge-Kutta		&\texttt{erk}\\
	\hline
	Crank-Nicolson			&\texttt{irk}\\
	\hline
	REXI						&\texttt{rexi}\\
	\hline
	Exponential Time Differencing Runge-Kutta	&\texttt{etdrk}\\
	\hline
\end{tabular}
\caption{\label{tbl:time_integrators_names}
	Description and string identifiers of different time integrators used throughout this work.
}
\end{table*}

We briefly introduce a string annotation to differentiate between the implementations of the time splitting methods.
Section \ref{sec:shallow_water_equations} introduced the different terms $L_g$, $L_c$, $N$ of the SWE PDE.
These terms will be denoted as \texttt{lg}, \texttt{lc}, \texttt{n}, respectively, and an overview is given in Table \ref{tbl:swe_terms}.
Abbreviations for different single-step time integrators used in this work are given in Table \ref{tbl:time_integrators_names}.
We use a concatenation of operators followed by the time stepping method to denote them.
E.g.\,using a Strang-splitting where the linear term is assumed to be $L_g$ and treated by a Crank-Nicolson scheme and the other non-linear terms assumed to be the Coriolis $L_c$ and non-linear terms $N$ and treated by an explicit \nth{2} order Runge-Kutta is denoted by \texttt{lg\_irk\_lc\_n\_erk\_ver0} where \texttt{ver0} denotes the version integrating the linear parts twice in each time step, see Sec.\,\ref{sec:strang_splitting}.

For ETD2RK, non-linearities are natively supported and no Strang-splitting is needed to gain a \nth{2} order accurate method.
This is denoted e.g.\,as \texttt{lg\_rexi\_lc\_n\_etdrk} which would treat the term $L_g$ as the only linear one with REXI and the terms $L_c + N$ as the non-linear ones within the ETD2RK time integration (see Sec.\,\ref{sec:etdnrk}).

\section{Results and discussions}
\label{sec:results}

Next we assess the performance of the previously discussed time stepping methods regarding various properties.
For the space-discretization, all studies are performed with SH using a T128 triangular truncation, see e.g.\,\cite{scott2004bob}, including a larger resolution in physical space to avoid anti-aliasing.

The barotropic instability benchmark\cite{galewsky2004initial} belongs to the most frequently used benchmarks to study new numerical schemes for the SWE on the sphere and uses parameters which generate simulation conditions related to the real atmosphere.
This benchmark is setup with a geostropic balanced initial condition in the northern hemisphere which is intended to remain by itself steady over the entire simulation runtime.
However, an intentional disturbance is added in form of a Gaussian bump which leads to well-defined slowly evolving vortices over several days, see \cite{galewsky2004initial} for a detailed description and a visualization.
This benchmark therefore allows us to study time integration methods over time integration ranges spanning multiple days and a non-trivial non-linear test case.
Additionally, we use zero-viscosity in the following benchmarks.
The error is computed using the $L_{\infty}$ norm on differences to the reference solution at $t=120h$.
This reference solution is based on a simulation using 4th order Runge-Kutta with a time step size of $\Delta t=15/8s$.

These studies are conducted on the Cheyenne supercomputer, equipped with 2.3-GHz Intel Xeon E5-2697V4 (Broadwell) processors,
18 cores per socket with dual-socket nodes.
Frequency scaling is activated per default on this system which is one reason for varying runtimes even for the same setup.
Therefore, we requested a fixed frequency of $2.3$ GHz.
With the library configuration playing an important role in the wallclock time performance, we'd like to briefly summarize this here:
\mbox{ICPC~18.0.1}, \mbox{SGI~MPT~2.15}, \mbox{FFTW~3.3.8} (\emph{precomputed \& optimized plans}), \mbox{SHTNS~3.1.0}, SWEET\cite{sweet_repository}~ver.\,2018-10-15.
The SWEET source code with the full implementation of this work is freely available at \cite{sweet_repository}.
The network also plays a crucial role for the parallel reduce in the REXI parts since communication of the full spatial data is required.
It is based on a Mellanox EDR InfiniBand high-speed interconnect with 18 nodes per switch node and a 25 GBps bidirectional link.
The coarser network topology is based on a partial 9D enhanced hypercube single-plane interconnect topology and is dimension-order routed.
Therefore, this yields efficient logarithmic $O(\log_2 n)$ complexity for reduce operations involving $n$ ranks.

For the REXI parameters, we start with $p_0 = 10.0$, $p_{\pm1} = \pm 30i$ and $N=128$ which we determined empirically (alternatives included $p_0 = 10.0 \pm1 5$ and $N=\{64, 128, 256\}$, $p_{\pm1}=30i \pm 10i$).
For REXI-based simulations using these parameters, time step sizes below $60\,\text{sec.}$ showed even an increase of the error, hence pointing out a missing consistency.
Similar properties were already identified in terms of the linear geostrophic balanced benchmark, see \cite{Schreiber2017SHREXI} based on the T-REXI coefficients.
Although the results were still competitive to other time integration methods, we gained a consistency for the considered time step sizes by scaling the radius for the contour linearly with the time step size, using circle contour radius of $30$ for a time step size of $480\,\text{sec.}$ as a baseline.
Additionally, we used a lower limit of the circle radius of $5.0$ and reduced the number of contour quadrature poles to $32$.
Finally, this led to the expected convergence also for small time step sizes with REXI.

\subsection{Performance aspects and parallelization in space}
\label{sec:parallelization_in_space}

\begin{figure}
	\begin{center}
		\includegraphics[width=\columnwidth]{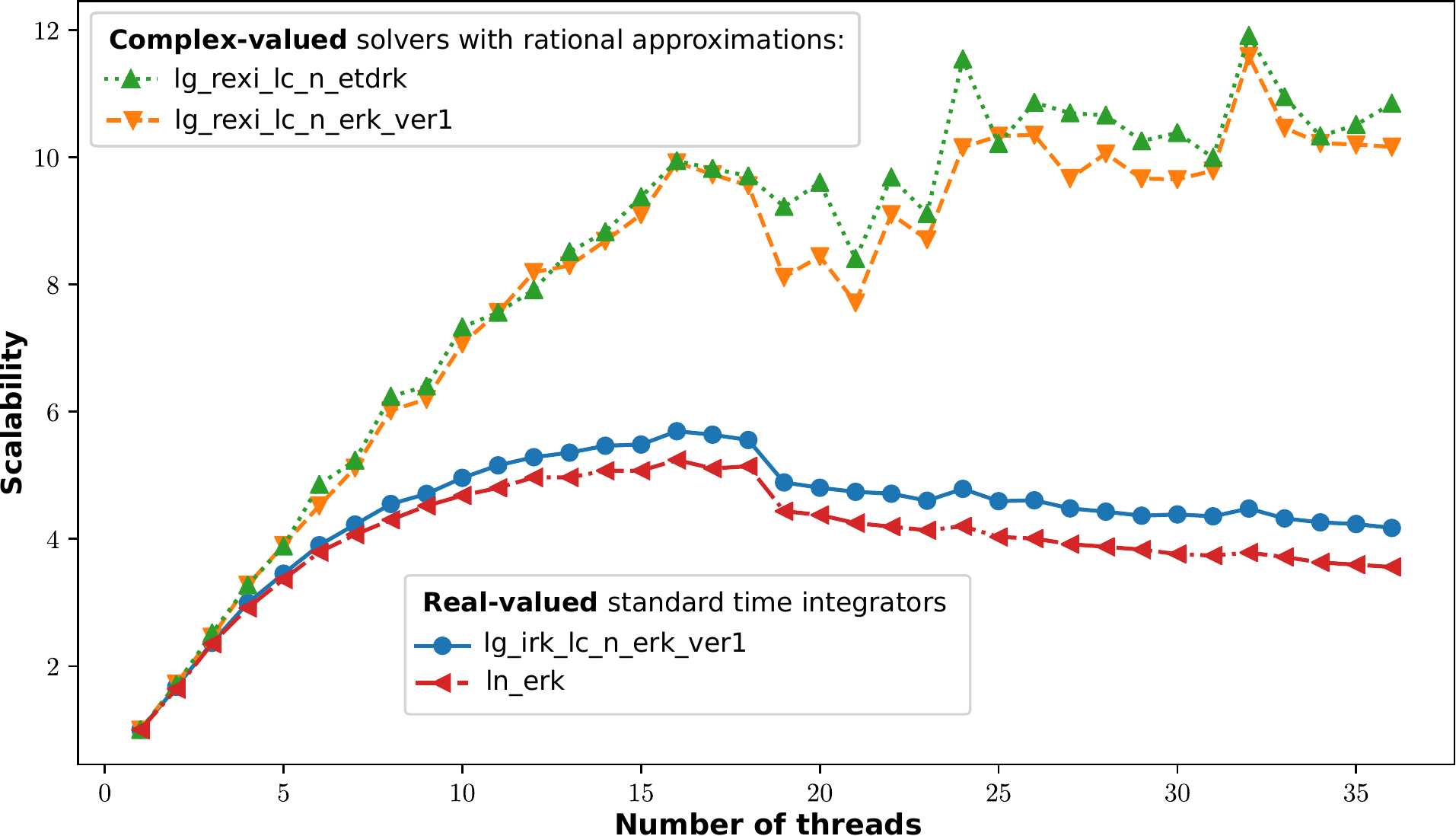}
	\end{center}
	\caption{\label{fig:output_threads_vs_scalability_edited}
		Scalability studies for
		explicit Runge-Kutta 2 (\emph{ln\_erk}),
		Crank-Nicolson \emph{lg\_irk\_lc\_n\_erk\_ver1}),
		Strang-split REXI (\emph{lg\_rexi\_lc\_n\_erk\_ver1}) and
		REXI-based ETD2RK (\emph{lg\_rexi\_lc\_n\_etdrk}) methods.
		We can observe a strongly improved scalability for the REXI solvers which can be related to the additional complex-valued workload.
		The scalability stagnates at about 18 cores which is related to the max.\,number of cores on the first socket.
	}
\end{figure}

We first study the performance of a parallelization in space to observe the scalability behaviour of standard time integration methods as well as the new ones with the results given in Figure \ref{fig:output_threads_vs_scalability_edited}.
The simulation executed 24 time steps and four different time integration methods were investigated:
explicit Runge-Kutta 2 (\emph{ln\_erk}, $0.354\,\text{sec.}$),
Crank-Nicolson \emph{lg\_irk\_lc\_n\_erk\_ver1}, $0.719\,\text{sec.}$.),
Strang-split REXI (\emph{lg\_rexi\_lc\_n\_erk\_ver1}, $31.66\,\text{sec.}$) and
REXI-based ETD2RK (\emph{lg\_rexi\_lc\_n\_etdrk}, $95.78\,\text{sec.}$).
The first two are classic ones and only involve real-valued computations.
The latter two methods involve computationally intensive REXI methods with 128 REXI terms.
The total wallclock time of the time stepping loop on one single core is also provided in the parentheses, showing that REXI with 128 terms takes about two order of magnitude longer to compute.
We like to emphasize that the parallelization-in-time of REXI across the 128 terms will be exploited later.

Regarding the scalability, the workload plays a crucial role for it:
although working with the same number of modes, REXI solvers involve complex-valued weights for all modes.
Therefore, this results in additional workload: for each of the 128 terms, twice as much bandwidth to load/store coefficients as well as additional complex-valued calculations even in physical space.
This additional workload of REXI-involving solvers eventually leads to the improved scalability 
which we will exploit later on.
No further significant scalability can be gained beyond $18$ cores which also represents the max.\,number of cores on the first socket which are used throughout the remainder of this work.
We emphasize that this scalability limit in the space domain depends on hardware architecture but also spatial resolution and model complexity: high resolution, full atmospheric dynamical cores with multiple layers (hence significantly increased per-cell workload in the vertical) will scale well beyond $18$ cores.

\subsection{Error vs. timestep size}

\begin{figure*}
	\begin{center}
		\includegraphics[width=0.7\textwidth]{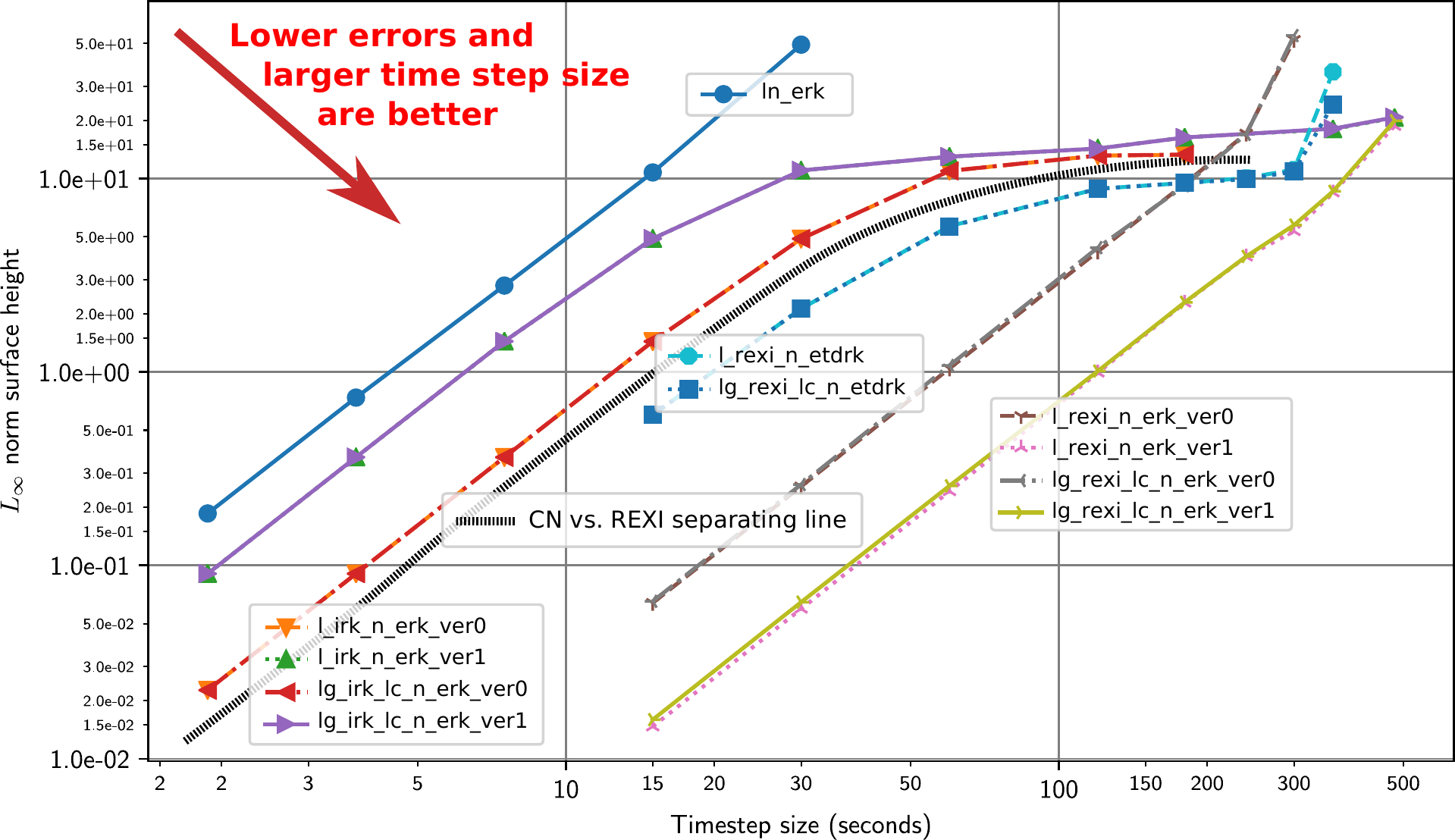}
	\end{center}
	\caption{\label{fig:barotropic_instability_h_dt}
		Error plots ($L_{\infty}$ norm on the reference solution) of the surface height based on the barotropic instability benchmark at simulation time $t=120h$ for varying time step sizes.
		We can observe a well-known significant time-step limitation of the explicit RK2 time stepping method and implicit one resulting in significantly larger stable time step sizes.
		Using ETD2RK scheme leads to reductions in the errors compared to the implicit time stepping methods for its stable time step regime.
		A Strang-split REXI leads to a max.\,stable time step size which similar to the implicit one.
		For smaller time step sizes, the best performing REXI method leads to significant reductions in the errors.
	}
\end{figure*}

Performance improvements are frequently first discussed in terms of time step sizes and we start by error vs.\,timestep size comparisons given in Fig.\,\ref{fig:barotropic_instability_h_dt}.
For purely explicit time integrators (\texttt{ln\_erk})we used time step sizes
$$\Delta t_{expl} \in \{ 15/8, 15/4, 15/2, 15, 30, 60 \},$$
for REXI-based ones we used
$$\Delta t_{rexi} \in \{15, 30, 60, 120, 180, 360, 480, 600\}$$
and
for implicit ones (\texttt{*\_irk\_*}) we tested with $\Delta t_{impl} \in t_{expl} \cup t_{rexi}$.
For sake of a better overview, we filtered results exceeding an $L_{\infty}$ error of $100m$ after 5 days of simulation.

We start with observations involving different types of time integrators.
Independent of the time stepping method, we can observe that \texttt{ver1} time split formulations allow larger stable time step sizes.
We account for this by time-integrating the non-linearities in each full time step with two time integrations, each one over half a time step size.
This can be also interpreted in a way that two non-linear integrations are executed in a row (at the end of one time step and the beginning of another time step), hence resulting in two time integrations of the non-linearities.

\begin{itemize}
	\item
	\textbf{Fully explicit Runge-Kutta \nth{2} order} (\texttt{*\_erk\_*}):\newline
	The Runge-Kutta time stepping method is strongly limited in their maximum time step size because of the fast gravity modes (see Sec.\,\ref{sec:timestepping_methods_for_c_and_w}) and stability issues of RK2 for oscillatory problems.

	\item
	\textbf{Implicit/explicit} (\texttt{*\_irk\_*}):\newline
	All \nth{2} order implicit time stepping methods show larger errors compared to REXI due to slowing down of fast propagating modes\cite{durran2010numerical}.
	Even more interesting once comparing to REXI, we observe its convergence in the time step size requires strongly reduced time step sizes compared to REXI.
	The maximum stable time step size of $\Delta t = 480s$ is similar to the one of REXI.
	In particular, both methods show a similar order of magnitude for the error.
	We account for the overlap of the \texttt{l\_irk\_*} and \texttt{lg\_irk\_*} lines by stiffness properties of the Coriolis effect allowing stable time step sizes beyond the one of the non-linearities.
	This Coriolis effect is treated in both cases with \nth{2} order accurate methods (RK and CN) and both of them are prone to similar dispersion errors \cite{durran2010numerical} in this case.
	Even if \texttt{*\_irk\_*ver1} allows to run stable for larger time step sizes, its overall error is higher compared to \texttt{ver0} formulations.

	\item
	\textbf{ETD2RK} (\texttt{*\_rexi\_*\_etdrk\_*}):\newline
	The ETD2RK scheme leads throughout all simulations to larger errors compared to the best Strang-split REXI.
	However, it also led to reduced errors for its stability region compared to the implicit time integrators.
	This is due to the accurate EI treatment of the linear terms.

	\item
	\textbf{Strang-split REXI} (\texttt{*\_rexi\_*\_erk\_*}):\newline
	Throughout all studies, a Strang-split REXI version provided the best error-vs.-timestep size results, resulting in a significant reduction of the errors.
\end{itemize}

\subsection{Error vs. wallclock time}

\begin{figure*}
	\begin{center}
		\includegraphics[width=0.7\textwidth]{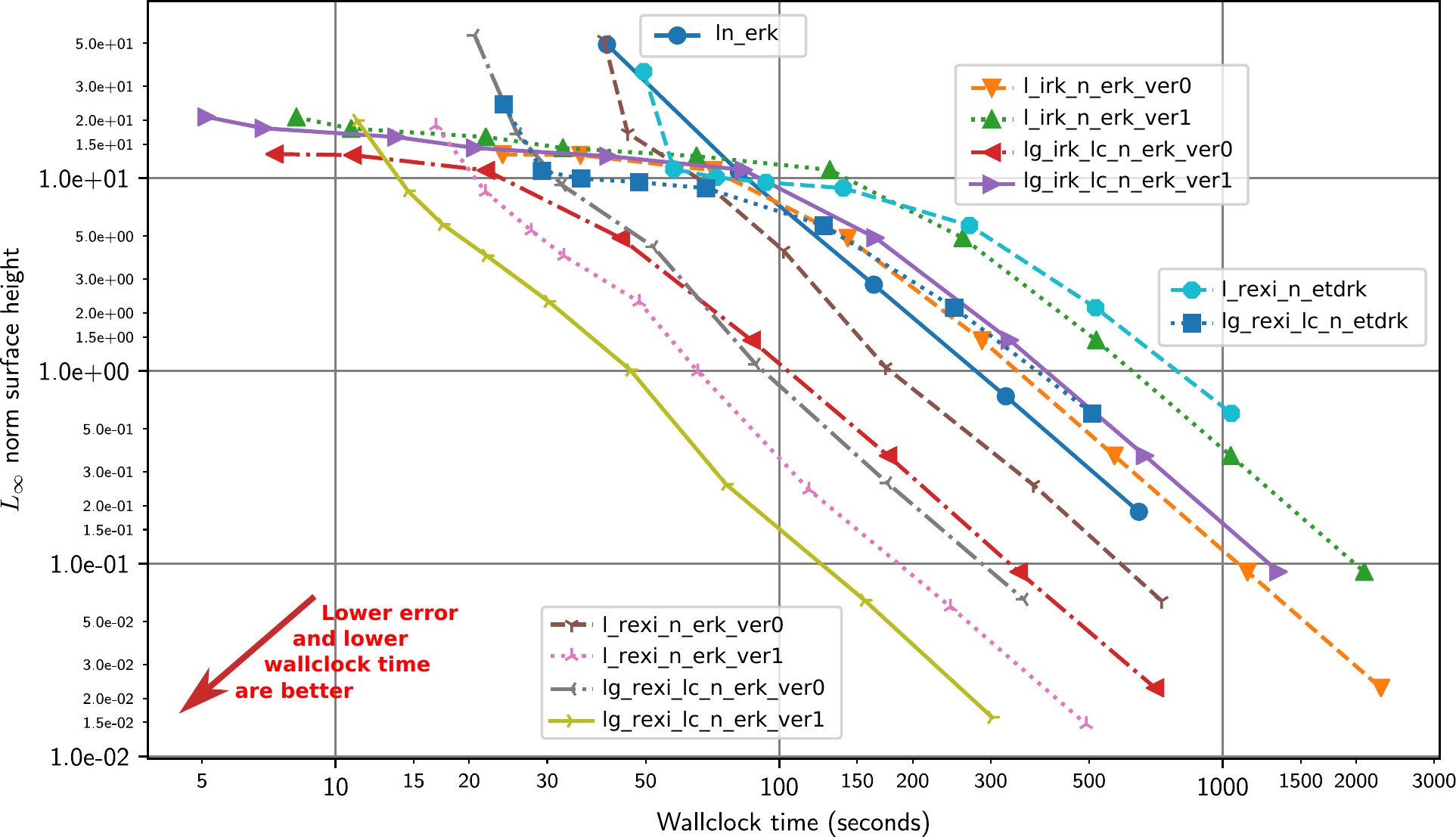}
	\end{center}
	\caption{
		\label{fig:barotropic_instability_wallclocktime}
		Wallclock time studies for the barotropic instability benchmark over a 5 day time integration interval.
		The explicit RK2 method (ln\_rk2) has in general a significantly larger wallclock-time to error ratio compared to the REXI and implicit time stepping solvers.
		Therefore, the implicit method is frequently preferred since this allows taking significantly larger time steps while still keeping the error on a moderate level.
		Regarding exponential integrators, none of the ETD2RK schemes are competitive to the best implicit Strang-split formulations.
		However, the Strang-split REXI method provides the best wallclock-time to error for more relaxed wallclock time restrictions and is Pareto optimal if a wallclock time of 100 seconds is affordable.
	}
\end{figure*}

Considering the errors with respect to the time-step size does not take into account the actual computational workload involved in each time step.
Therefore, we additionally conducted performance studies on a supercomputer comparing the wallclock time vs.\, the errors with the results given in Fig.\,\ref{fig:barotropic_instability_wallclocktime}.

We used one MPI rank per available REXI term, hence using for e.g. $N=128$ REXI terms $64$ compute nodes and each rank's threads are pinned to the cores of each rank-dedicated socket.
For the parallelization-in-time we solely use MPI for the broadcast and reduction operators over the rational approximations with each rank's master thread.

\begin{itemize}
	\item
	\textbf{Fully explicit Runge-Kutta \nth{2} order} (\texttt{*\_erk\_*}):
	Despite the computational simplicity of RK2, the small time-step sizes which are enforced by the RK2 method also lead to larger wall clock times, making it not competitive to e.g.\,the implicit time stepping method which is, for similar errors, about $3 \times$ faster.

	\item
	\textbf{Implicit/explicit} (\texttt{*\_irk\_*}):
	The implicit time stepper (Crank-Nicolson) results in the largest time step sizes and, because of its computational simplicity, also to the lowest wallclock time for a stable time integration.
	However, we can also observe a significant increase in errors for lower wallclock time, hence larger time step sizes.
	Again, this is due to the fast propagating modes which are increasingly slowed down with increasing time step sizes (see \cite{durran2010numerical}).
	We like to mention here that we exploited the possibility of solving for the system of equations directly in spectral space, but it might be more time consuming using non-spectral methods.

	\item
	\textbf{ETD2RK} (\texttt{*\_rexi\_*\_etdrk\_*}):
	Taking only the time step size into account (see previous section), this scheme provided partly improved results.
	However, the wallclock time reveals that it is not competitive to the other methods:
	in all our simulations, the ETD2RK scheme never provided improved wallclock time vs.\,error compared to the best one of the other time integration methods.

	\item
	\textbf{Strang-split REXI} (\texttt{*\_rexi\_*\_erk\_*}):
	We put our focus on the best REXI time stepping method \linebreak[4](\texttt{lg\_rexi\_lc\_n\_erk\_ver1}):
	Here, compared to all other methods the errors are significantly reduced for the stable timestep size regime.
	Again, in contrast to Crank-Nicolson, REXI does not suffer from slowing down the fast propagating waves.
	Because of the additional computations and parallelization overheads involved in REXI, the wallclock time is larger if compared to the implicit time stepping method and the same time step size.

\end{itemize}

Next, we compare the best REXI method (lg\_rexi\_lc\_n\_erk\_ver1) with the best implicit method (lg\_irk\_lc\_n\_erk\_ver0) from two different perspectives:
The first one is a requirement on certain accuracy of the simulations.
With such a requirement of e.g. $L_\infty(h-h_{ref}) \leq 0.1$, we can observe an almost $3 \times$ improved wallclock time of the best REXI method compared to the best implicit time stepping method.
Similarly, with a second perspective, one can request the simulation to be finished within a certain time frame.
With e.g.\,$100\,\text{sec.}$, we can observe an over $6\times$ reduced error for REXI compared to the next best method.
However, we should also point out a case where it is desirable to get results as fast as possible. In this case, the implicit time integration method provides the fastest results in under $10\,\text{seconds}$.

We also compare the treatment of the Coriolis effect ($L_c$) as part of the linear or non-linear time stepping method:
Whereas the results in Fig.\,\ref{fig:barotropic_instability_h_dt} indicated superior properties by treating the Coriolis effect as part of REXI, the wallclock time results in Fig.\,\ref{fig:barotropic_instability_wallclocktime} also incorporate the real computational costs.
Here it is important to mention that treating the Coriolis effect as part of the linear part required additional computations such as solving a penta-diagonal matrix instead of a diagonal-only, see \cite{Schreiber2017SHREXI}.
We can conclude that for this particular benchmark and using the errors on the height field, treating the Coriolis effect as part of the non-linearities provides the best REXI-based wallclock-time to error results.

\subsection{REXI parallelization overheads}

\begin{figure*}
	\begin{center}
		\includegraphics[width=0.76\textwidth]{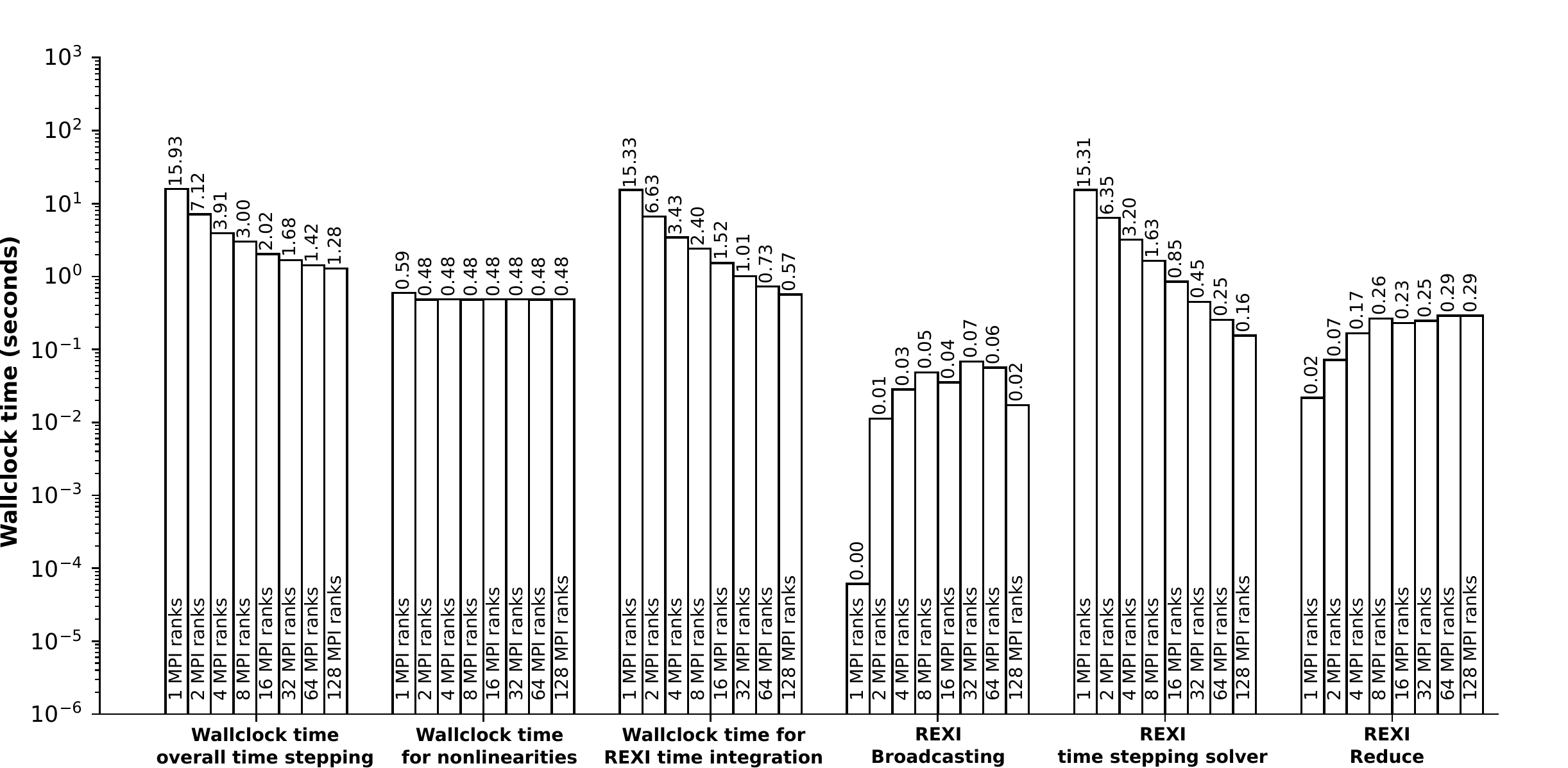}
	\end{center}
	\caption{\label{fig:output_rexi_performance_breakdown_postprocessed}
	Parallelization overheads of the REXI time stepping method (\texttt{lg\_rexi\_lc\_n\_erk\_ver1}) over 100 time steps.
	Wallclock time is plotted with log-scaling along y-axis.
	Each block shows the wallclock times for particular parts of the time integration over an increasing number of MPI ranks.
	We can observe that even assigning only one single REXI term to one rank, hence the extreme scale for this parallel-in-time method, the overall wallclock time is not dominated by the collective operations.
	}
\end{figure*}

We close the results with a detailed study on a wallclock time performance breakdown of the REXI time stepping method.
We timed different operations using the \texttt{lg\_rexi\_lc\_n\_erk\_ver1} time stepping method with a stable time step size of $\Delta t = 300s$.
Only $100$ time steps were executed, hence only over a very short time integration range to reduce the likelihood of network package collisions on the shared network during one execution.
Even though, the generated results on Cheyenne were still prone to a large variance, mitigated using 10 ensemble runs.
For our plots, we used the minimum over individual benchmark values of all ensemble runs.
Here, we assume that these values represent the ones which would be generated on a system with fully exclusive resources, e.g.\,without any network congestion.
The results are given in Fig.\,\ref{fig:output_rexi_performance_breakdown_postprocessed}.

The first block (\emph{Wallclock time overall time stepping}) shows the total wall clock time to execute the time integration.
This time robustly decreases for increasing number of MPI ranks.

In the second block (\emph{Wallclock time for nonlinearities}) we can observe a constant wallclock time to evaluate the non-linear parts of REXI.
These parts don't involve any communication and remain constant throughout the full REXI scaling studies.

The third block (\emph{Wallclock time for REXI time integration}) depicts the wallclock times only for the REXI-related parts.
Overall, this shows a robust performance improvement for increasing the number of MPI ranks to parallelize over.
The last bar shows a total runtime of $0.57\,\text{sec.}$, hence a very similar wallclock time compared to the one evaluating the non-linearities.

The fourth block (\emph{REXI Broadcasting}) depicts the timings for the broadcast operation.
This was measured on the \nth{2} rank to avoid any potential fire-and-forget behaviour of the MPI implementation on the first rank.
Besides minor fluctuations \footnote{Still induced by natural effects on shared network resources.} we can observe a very low broadcast overhead of almost constant runtime throughout all test cases.
Also, the broadcast operation is significantly faster than the parallel reduce.

The fifth  block (\emph{REXI time stepping solver}) shows the wallclock time required to solve for all the REXI terms.
Here, we can observe a robust scaling and we like to point out that with 128 terms in the REXI sum, the maximum scalability for the parts related to the REXI parallelization is also reached with 128 MPI ranks.

The last block (\emph{REXI Reduce}) shows overheads of the reduce operation which are increasing with the number of ranks.
However, we can also observe an $O(\log(N))$ communication complexity of the underlying hypercube interconnect.
Therefore, we don't expect that the communication overheads will further significantly dominate even for a larger number of REXI terms.

We discuss these REXI timings with the best \nth{2} order Strang-split Crank-Nicolson method lg\_irk\_lc\_n\_erk.
This method requires $0.572\,\text{sec}$ for the similar time integration range whereas REXI on its full scaleout requires $1.28\,\text{sec.}$
Therefore, REXI in this formulation is $2.2 \times$ more expensive than the Strang-split CN method.
However, we like to point out that REXI also provides higher accuracy.

Finally, we briefly discuss the theoretical scalability of REXI-based time stepping methods, 
Since our studies are conducted with a fixed workload, this makes Amdahls' law applicable.
Here, the scalability limitations are induced by the sequential parts which are in our case the treatment of the non-linearities and the pre- and post-processing of REXI.
These sequential parts obviously dominate, hence also limit, the strong scalability of REXI.
Again, we would like to emphasize that PinT scaling does not follow the standard expectations from spatial scalability, but that every additional benefit counts, which is the case here even for the maximum theoretical scalability of $128$ ranks.

\section{Summary and discussion}

The motivation of this paper is to improve weather and climate simulations.
In that context, we have compared the competitiveness of parallel-in-time rational approximations of exponential integrators (REXI) methods, based on Cauchy contour integrals, to more conventional time integration methods.
We conducted studies with the often-used barotropic instability test case for the non-linear shallow-water equation on the rotating sphere to assess the performance of each time integration method.

Previous work showed that REXI allows arbitrarily long time-steps for the linearized equations.
In the fully nonlinear case, we use REXI to integrate the linear parts of the shallow-water equations either using a time splitting approach or by using the ETDnRK exponential integrator method.
Additional challenges are posed by using the Cauchy contour integral to derive REXI coefficients, including the requirement of large time step sizes and avoiding numerical cancellation effects.
We circumvented these issues by modifying the integration contour and by focusing on physically relevant linear oscillatory and diffusive stiff problems.

We first compared the errors for each time integration method and the time step size, revealing the superior properties of REXI in almost all cases.
However, this approach ignores the actual computational performance of each method.
Therefore, we compared integration rate vs.\,error on NCAR's Cheyenne supercomputer.
On the one hand, focusing solely on fastest time-to-solution for a stable time integration method, the REXI-based method was not able to compete with an implicit time stepping method, due its additional computational costs.
On the other hand, we also studied REXI's stable time step regime with a comparison to the best implicit time integration method.
Demanding for a low error requirement, REXI showed a close-to threefold reduction in wallclock time.
Alternatively, requirements on finishing the simulation run within a particular wallclock time frame led to an over $6\times$ reduced error using REXI.

Finally, we conducted a detailed performance analysis of REXI to gain insight into its parallelization overhead.
As REXI requires to communicate the entire full-dimensional volume via MPI broadcast and reduce operations, it might be assumed to significantly limit its scalability.
Even so, our benchmarks reveal that the communication overhead along the time dimension is not dominating the overall simulation runtime to
solve the non-linear equations on the Cheyenne supercomputer.
An optimization of the reduce operation is expected to lead to valuable improvements of the parallel REXI time integration method.

Our results strongly suggest that REXI is a good candidate time integration method for climate and weather simulation code, because of its ability to overcome the oscillatory linear stiffness in the SWEs without introducing time integration errors in the linear parts.
	Regarding the parallelization, such a spatial workload can be treated orthogonal to the REXI time-parallelism.
	Therefore, we expect additional scalability once the spatial workload is increased.
	This is the case e.g. for full dynamical cores for climate and weather simulations requiring computations on larger spatial workloads (additional vertical discretization and possibly higher horizontal resolution).
	Regarding the interconnect and considering the return of investments on supercomputers, we envision heterogeneous interconnects and corresponding scheduler support to exploit the different communication properties in time and space.

For future work, we see research on the treatment of the non-linear terms with, for example a semi-Lagrangian method, which might be expected to lead to further significant reductions in errors compared to implicit schemes.
Ways to further treat the remaining non-linearities in time in the context of simulations which are relevant for climate and weather simulations include e.g.\,ML-SDC method\cite{Francois2018MLSDC} as well as the PFASST method\cite{emmett2012toward}.

Since the number of REXI terms directly relates to the parallel overheads and, even worse, linearly to the total computational costs, we see a reduction of the number of these terms as an opportunity to further reduce these computational costs, hence making REXI-based methods more practical in
operational atmospheric models.

\section*{Acknowledgements}

Martin Schreiber received funding from NCAR for a research stay in summer 2017 at the Mesa Labs.
We'd like to acknowledge Cheyenne \cite{cheyenne} supercomputer used to assess the wallclock time performance and also the ``CoolMUC Cluster at LRZ, Germany'', the ``MAC cluster at Technical University of Munich'' as well as ``The Applied Mathematics Computational Laboratory of the Institute of Mathematics and Statistics, University of São Paulo'' for providing computational resources for all other studies.

\bibliographystyle{elsarticle-num-names}
\bibliography{bibliography}

\end{document}